\documentclass[12pt]{amsart}

\usepackage{graphicx}

\numberwithin{equation}{section}
\usepackage[T1]{fontenc}

\graphicspath{ {images/} }
\bibliography{sources}

\newcommand{\ita}{\textit}

\newcommand{\N}{\mathbb{N}}
\newcommand{\Z}{\mathbb{Z}}
\newcommand{\C}{\mathbb{C}}
\newcommand{\Q}{\mathbb{Q}}

\newcommand{\R}{\mathbb{R}}

\newcommand{\e}{\epsilon}

\newcommand{\sub}{\subseteq}

\newcommand{\dist}{\mathrm{dist}}

\newcommand{\lb}{\left\langle}
\newcommand{\rb}{\right\rangle}

\newcommand{\ga}{\gamma}

\newcommand{\la}{\lambda}

\newcommand{\bpf}{\begin{proof}}
\newcommand{\epf}{\end{proof}}
\newcommand{\bal}{\begin{align*}}
\newcommand{\eal}{\end{align*}}
\newcommand{\baln}{\begin{align}}
\newcommand{\ealn}{\end{align}}
\newcommand{\beqn}{\begin{equation}}
\newcommand{\eqen}{\end{equation}}

\newcommand{\bnum}{\begin{enumerate}}
\newcommand{\enum}{\end{enumerate}}
\newcommand{\bpm}{\begin{pmatrix}}
\newcommand{\epm}{\end{pmatrix}}
\newcommand{\T}{\theta}

\newcommand{\A}{\alpha}

\newcommand{\B}{\beta}
\newcommand{\D}{\Delta}
\newcommand{\de}{\delta}

\newcommand{\Sym}{\mathrm{Sym}}
\newcommand{\s}{\sigma}

\newcommand{\bfr}{\begin{flushright}}
\newcommand{\efr}{\end{flushright}}

\newcommand{\Ga}{\Gamma}

\newcommand{\sm}{\setminus}

\newcommand{\set}[1]{\left\{ #1 \right\}}

\newcommand{\bitem}{\begin{itemize}}
\newcommand{\eitem}{\end{itemize}}

\newcommand{\paren}[1]{\left( #1 \right)}

\newcommand{\cO}{\mathcal{O}}

\newcommand{\inv}{^{-1}}

\newcommand{\Mod}[1]{\ (\mathrm{mod}\ #1)}
\newcommand{\W}{\Omega}

\newcommand{\leg}[2]{\left(\frac{#1}{#2}\right)}

\newcommand{\cD}{\mathcal D}

\newcommand{\cl}{\overline}

\newcommand{\btikz}{\begin{tikzcd}}
\newcommand{\etikz}{\end{tikzcd}}

\newcommand{\rabs}[1]{\left|#1\right|}

\DeclareMathOperator{\PGL}{PGL}

\renewcommand{\phi}{\varphi}

\newcommand{\Zpos}{\mathbb{Z}_{>0}}

\renewcommand{\phi}{\varphi}

\renewcommand{\phi}{\varphi}

\newcommand{\bcases}{\begin{cases}}
\newcommand{\ecases}{\end{cases}}

\DeclareMathOperator{\Gal}{Gal}

\usepackage{relsize}
\usepackage{mathrsfs}

\newcommand{\n}[1]{\left\| #1 \right\|}

\newcommand{\cP}{\mathcal P}
\newcommand{\Pos}{\mathrm{Pos}}
\newcommand{\cM}{\mathcal{M}}
\newcommand{\Mat}{\mathrm{Mat}}
\newcommand{\den}{\mathrm{den}}
\newcommand{\cS}{\mathcal S}
\newcommand{\cL}{\mathcal L}
\newcommand{\cQ}{\mathcal Q}
\newcommand{\cH}{\mathcal H}
\DeclareMathOperator{\PO}{PO}
\DeclareMathOperator{\SL}{SL}
\DeclareMathOperator{\rk}{rk}
\newcommand{\bs}{\backslash}
\DeclareMathOperator{\sgn}{sgn}

\RequirePackage[colorlinks, linkcolor=blue, citecolor=blue, urlcolor=blue]{hyperref}

\RequirePackage[margin=1in]{geometry}

\theoremstyle{plain}

\newtheorem{theorem}{Theorem}[section]
\newtheorem{lemma}[theorem]{Lemma}
\newtheorem{proposition}[theorem]{Proposition}

\newtheorem*{remark*}{Remark}
\newtheorem*{remarks*}{Remarks}

\usepackage{mathtools}

\numberwithin{equation}{section}

\begin{document}

\title{Explicit subconvexity savings for sup-norms of cusp forms on $\mathrm{PGL}_n(\R)$}

\author{Nate Gillman}
\address{Department of Mathematics \& Computer Science, Wesleyan University, Middletown, CT 06457, U.S.A.}
\email{ngillman@wesleyan.edu}

\date {June 17, 2019}

\begin{abstract}
Blomer and Maga \cite{BlomerMaga} recently proved that, if $F$ is an $L^2$-normalized Hecke-Maass cusp form for $\mathrm{SL}_n(\mathbb Z)$, and $\Omega$ is a compact subset of $\mathrm{PGL}_n(\mathbb R)/\mathrm{PO}_n(\mathbb R)$, then we have $\|F|_\Omega\|_\infty\ll_\Omega\lambda_F^{n(n-1)/8-\delta_n}$ for some $\delta_n>0$, where $\lambda_F$ is the Laplacian eigenvalue of $F$.
In the present paper, we prove an explicit version of their result.
\end{abstract}
\maketitle


\section{Introduction and Statement of Results}

An automorphic form $F$ is defined on a quotient $\Ga\bs X$ of a Riemannian symmetric space by a discrete subgroup of its isometries.
A fundamental property of an automorphic form is its size, and in particular the distribution of its mass.
One measure of equidistribution is a bound of some $L^p$-norm of $F$, an especially important case being $p=\infty$.
As an automorphic form is an eigenfunction of the Laplacian, of particular interest is bounding a given automorphic form in terms of its Laplacian eigenvalue $\la_F$.
In 2004, Sarnak \cite{Sarnak1} proved that, if $X$ is a compact locally symmetric space and $\mathscr D(X)$ is the algebra of differential operators invariant under the Riemannian isometry group of $X$,
then an $L^2$-normalized joint eigenfunction $F$ of $\mathscr D(X)$ satisfies the following bound,
\begin{equation}\label{eq:Sarnakbound}
\n F_\infty\ll\la_F^{(\dim X-\rk X)/4}.
\end{equation}
This result, which was proved using purely analytic arguments, is often referred to as the \ita{convexity bound}, and it is known that the exponent is sharp in general.

Here we are interested in arithmetic situations.
Many classical examples of Riemannian locally symmetric spaces enjoy additional symmetries given by the Hecke operators, a commutative family of ``averaging'' operators that play an important role in the theory of modular and automorphic forms; see for example \cite{Ono}.
In these situations, automorphic forms on $X$ are also joint eigenfunctions of the Hecke algebra.
In light of this additional layer of symmetry, it is reasonable to expect some power saving in (\ref{eq:Sarnakbound}) when we restrict $F$ to compact subsets of $X$. 
Such a restriction is necessary in order to avoid large growth at cuspidal regions, see for example \cite{BrumleyTemplier}.
This is often referred to as the \ita{subconvexity conjecture for sup-norms of cusp forms.}

The first discovery of subconvexity is due to Iwaniec and Sarnak \cite{IwaniecSarnak} in 1995.
They demonstrated a saving of $1/24$ for automorphic forms on the hyperbolic plane $\cH_2$.
Since then, much work has been done in this area, but only recently has any power-saving been discovered for higher rank spaces: 
in 2014, Blomer and Pohl \cite{BlomerPohl} proved subconvexity for Hecke-Maass cusp forms on the Siegel modular space of rank $2$; see also \cite{HRR} and \cite{BlomerMaga2}.
Additionally, a preprint of Marshall \cite{Marshall} demonstrates a power saving for a wide class of semi-simple groups.

In 2016, Blomer and Maga \cite{BlomerMaga} proved subconvexity for Hecke-Maass cusp forms on $\SL_n(\Z)$, for all $n$.
They provided a proof of some power saving without explicating it. 
Until this paper, no explicit power saving has been given for the cases $n\geq 3$ in this general setting, which is our main result.

\begin{theorem}\label{thm:1}
Let $n\geq 2$.
Let $F$ be an $L^2$-normalized Hecke-Maass cusp form on $\SL_n(\Z)$, and let $\W$ be a fixed compact subset of $\PGL_n(\R)/\PO_n(\R)$.
Then,
\begin{equation}\label{eq:mainthmFbound}
\n{F|_\W}_\infty\ll_{\W,\e}\la_F^{\frac{n(n-1)}8\,-\,\de_n\,+\,\e},
\end{equation}
where $\de_n \gg n^{-cn^6}$ is given explicitly by (\ref{eq:delta_ndefn}).
\end{theorem}

In Table \ref{table:1} we provide numerical values for the first few $\de_n$.
Our proof gives an exact formula for these, but does not optimize the value.
\begin{table}[h] 
\begin{center}
\begin{tabular}{|c|c|c|c|c|c|c|}\hline 
    $n$ & $\de_n\approx$ \\\hline 
    $2$ & $2.74\cdot10^{-143}$\\\hline
    $3$ & $6.44\cdot10^{-976}$\\\hline 
    $4$ & $2.29\cdot10^{-3951}$\\\hline
    $5$ & $2.39\cdot10^{-12273}$\\\hline
    $6$ & $4.71\cdot10^{-32175}$\\\hline
	$7$ & $9.58\cdot10^{-74679}$\\\hline
	$8$ & $9.28\cdot10^{-157867}$\\\hline
\end{tabular}
\end{center}
\caption{Approximate values given by Theorem \ref{thm:1}}
\label{table:1}
\end{table}
It is humorous to compare our colossally small value of $\de_2$ against Iwaniec and Sarnak's \cite{IwaniecSarnak} breakthrough result of $\de_2=1/24$.

Using a different method, we also prove the following better bound in the case $n=3$.

\begin{theorem}\label{thm:2}
$\de_3=1/812$ is also suitable in (\ref{eq:mainthmFbound}).
\end{theorem}

This constant is optimized within the framework of our argument.
We should note that Holowinski, Ricotta, and Royer \cite{HRR} proved a result analogous to Theorem \ref{thm:2}, but in a more restricted setting.
Specifically, they proved $\de_3=1/76$ suffices, provided that the Hecke-Maass cusp forms have one Langlands parameter which is uniformly bounded.

\begin{remark*}
A slight modification of our argument gives an identical bound for $\n{F|_\W}_\infty$ in terms of spectral parameters; see the introduction to \cite{BlomerMaga}.
\end{remark*}

\begin{remark*}
Another small modification of our argument gives a nearly identical bound for Hecke-Maass cusp forms on a given congruence subgroup $\Ga_0(N)$; again, see \cite{BlomerMaga}.
\end{remark*}

Our paper is organized as follows.
In Section \ref{sec:preliminaries}, we recall a result from \cite{BlomerMaga} and explain a matrix-counting problem whose solution yields the proofs of Theorems \ref{thm:1} and \ref{thm:2}.
In Section \ref{sec:lemmas}, we prove technical Lemmas \ref{lem:Galois_bounding}, \ref{lem:peter5a}, \ref{lem:peter5b}, and \ref{lem:peter8}, involving diophantine analysis-style bounds over algebraic number fields, as well as Lemma \ref{lem:primesinintervals}, which provides an estimate on the quantity of primes in relevant dyadic intervals.
These results provide explicit bounds needed in Section \ref{sec:recursion}, where we prove Proposition \ref{prop:recursion}, which yields a good estimate for the matrix-counting problem.
We apply this in Section \ref{sec:proofofthm1} to prove Theorem \ref{thm:1}, and in Section \ref{sec:proofofthm2} we provide a proof of Theorem \ref{thm:2} using more elementary means.

\section{The counting problem} \label{sec:preliminaries}

For any $\ga\in\Mat_n(\Z)$, denote by $\D_j:=\D_j(\ga)$ the $j$'th determinantal divisor, i.e., the greatest common divisor of all $j\times j$ minors.
For any set of positive-definite matrices $\cM\sub\Mat_n(\R)$, for $Q\in\cM$, any $a,b\in\N$ and $M>0$ we define the following collection of integral matrices,
\begin{equation}\label{eq:Sdefn}
S(Q,a,b,M):=\set{\ga\in\Mat_n(\Z):
\ga^\top Q\ga-(ab^{n-1})^{2/n}Q\ll_\cM(ab^{n-1})^{(2-M)/n},
\D_1=1,
\D_2=b}.
\end{equation}
Here and throughout, we take estimates on matrices, such as the one above, to be entrywise.
Also, we formally allow $M=\infty$ to signify zero error.

By \cite[(2.8)]{BlomerMaga}, we have the following estimate.

\begin{proposition}\label{prop:blackbox}
Fix $n\geq 2$.
Let $\cL>5$ and $M>0$, and let $\cP$ be a set of primes in $[\cL,2\cL]$.
For $g\in\PGL_n(\R)$, define $Q:=(\det g)^2g^{-\top}g\inv\in\Sym_n(\R)$.
Let $F$ be an $L^2$-normalized Hecke-Maass cusp form on $\SL_n(\Z)$, and denote by $\la_F$ the corresponding Laplacian eigenvalue.
Then,
\begin{equation}\label{eq:blackbox}
\rabs{F(g)}^2\ll \la_F^{\frac{n(n-1)}4}
\paren{
\frac{1}{\rabs\cP}
+\la_F^{-1/4}\cL^{n^3+M/2}
+\sum_{\nu=1}^n\frac{1}{\rabs\cP^2}\sum_{p,q\in\cP}\frac{\#\cS(Q,q^\nu,p^\nu,M)}{\cL^{\nu(n-1)}}
}.
\end{equation}
\end{proposition}

The next two sections are devoted to bounding the cardinality of $\cS(Q,q^\nu,p^\nu,M)$, which is the matrix-counting problem discussed earlier.
Throughout, our argument uses that $g$ is in a fixed compact domain $\W$ of $\PGL_n(\R)/\PO_n(\R)$, so that for instance the implied constants in (\ref{eq:Sdefn}) and (\ref{eq:blackbox}) depend on $\W$ but not on $Q$.
Additionally, we take all implied constants to hold for sufficiently large $\cL$; this is acceptable because in our application of Proposition \ref{eq:blackbox}, we will take $\cL$ to be an increasing function of $\la_F$, and it is known that there are only finitely many Hecke-Maass cusp forms with bounded Laplace eigenvalue.
We also allow all implied constants in the first five sections to depend on $n$; we will indicate this dependence explicitly when it makes the argument more clear.


\section{Technical lemmas}\label{sec:lemmas}

In this section, towards estimating the cardinality of (\ref{eq:Sdefn}), we first prove four diophantine analysis-style lemmas.
Lemma \ref{lem:Galois_bounding} provides a Galois-theoretic framework for keeping track of error estimates.
We will repeatedly use this lemma when we prove Lemmas \ref{lem:peter5a}, \ref{lem:peter5b} and \ref{lem:peter8}, which are based on \cite[Lemmas 5(a), 5(b), 8]{BlomerMaga}. 
Our addition to these results is the incorporation of a scheme which makes the bounds proved in that paper's lemmas effective.

Let $K\sub\C$ be a number field and $\cl K\sub\C$ its Galois closure, and suppose $A>1$.
For a fixed number $\A\geq 1$, we say an element of $K$ is $\A$-well-balanced, or that it has well-balanced constant $\A$, if it can be written as a fraction $a/b$ with $a,b\in\cO_K$ and either $a=0$ and $b=1$, or else for each $\s\in\Gal(\cl K/\Q)$, we have
\[A^{-\A}\leq|\s a|,|\s b|\leq A^\A.\]

\begin{lemma}\label{lem:Galois_bounding}
Fix a number field $K$ and define $d_K:=\deg(\cl K/\Q)$.
If $a/b$ is $\A$-well-balanced and $c/d$ is $\B$-well-balanced, then:
\bnum
	\item The negation $-a/b$ has well-balanced constant $\A$.
	\item If $a\neq 0$, then the reciprocal $b/a$ has well-balanced constant $\A$.
	\item The product $ac/bd$ has well-balanced constant $\A+\B$.
	\item The sum $a/b+c/d$ has well-balanced constant $(\A+\B+1)d_K$.
\enum
Furthermore, the sum of $k$ elements of $K$, each with well-balanced constant $\A$, has the following well-balanced constant,
	\begin{equation}\label{eq:SkA}
	S_{k,d_K}(\A)=\begin{cases}
	k(\A+1)-1&d_K=1\\
	\A d_K^{k-1}+d_K(\A+1)\frac{d_K^{k-1}-1}{d_K-1}&d_K>1.
	\end{cases}
	\end{equation}
\end{lemma}
\bpf
The first three points are clear.
If $ad+bc=0$ then the fourth claim is obvious, so assuming $ad+bc\neq 0$, we can estimate
\[\prod_{\s\in\Gal(\cl K/\Q)}\rabs{\s(ad+bc)}
=\rabs{\mathcal N(ad+bc)}\geq 1.\]
It is clear that $\s(ad+bc)\leq A^{\A+\B+1}$, so for any $\s_0\in\Gal(\cl K/\Q)$, we have
\[\rabs{\s_0(ad+bd)}
=\rabs{\frac{\mathcal N(ad+bc)}{\prod_{\s\in\Gal(\cl K/\Q)\sm\set{\s_0}}\s(ad+bc)}}
\geq \frac{1}{A^{d_K(\A+\B+1)}},\]
which proves the fourth statement.
Therefore, when we sum $k$ terms, each with well-balanced constant $\A$, the well-balanced constant $S_{k,d_K}(\A)$ is given by the following linear recurrence,
\[S_{k,d_K}(\A)=d_K\cdot S_{k-1,d_K}(\A)+d_K(\A+1),\qquad S_{1,d_K}(\A)=\A,\]
which has closed form given in (\ref{eq:SkA}).
\epf

\begin{lemma}\label{lem:peter5a}
Let $m,r\in\N$ and $A\geq 2$.
Let $K\sub\R$ be a real number field and $\cl K\sub\C$ its Galois closure.
For $1\leq j\leq r$ let $b_j=(b_{1j},\dots,b_{mj})^\top\in\R^m$ and assume that all $b_{ij}$ are in the ring of integers $\cO_K$, and satisfy $b_{ij}=0$ or
\begin{equation}\label{eq:WBCdefn}
A\inv\leq|\s(b_{ij})|\leq A
\end{equation}
for all $\s\in\Gal(\cl K/\Q)$.
Let $H:=\cap_j b_j^\perp$.
Then for every $v\in\R^m,$ we have
\begin{equation}\label{eq:distvH}
\mathrm{dist}(v,H)\leq A^{\T_1}\max_j|\lb v,b_j\rb|,
\end{equation}
where $\T_1:=\T_1(m,r,d_K)$ is defined in (\ref{eq:theta1}) below.
\end{lemma}
\bpf
Take a maximal independent subset $\set{u_1^\top,\dots,u_{m'}^\top}\sub\set{b_1^\top,\dots,b_r^\top}$ (i.e., $\dim H=m-m'$.)
Then $u_1,\dots,u_{m'}$ is a basis in $H^\top$.
By the Gram-Schmidt procedure, we obtain inductively an orthogonal basis $u_j':=u_j-\sum_{i=1}^{j-1}\frac{\lb u_j,u_i'\rb}{\lb u_i',u_i'\rb}u_i'$ with entries in $K$.
The distance (\ref{eq:distvH}) is the following quantity,
\begin{align}\label{eq:distcomp}
\dist(v,H)
=\n{\mathrm{proj}_{H^\perp}v}
=\n{\sum_{j=1}^{m'}\frac{\lb v,u_j'\rb}{\lb u_j',u_j'\rb}u_j'}.
\end{align}
Say $u_j'$ has well-balanced constant $w_j:=w_j(m,r,d_K)$.
By Lemma \ref{lem:Galois_bounding}, $\frac{\lb u_j,u_i'\rb}{\lb u_i',u_i'\rb}u_i'$ has well-balanced constant $S_{m,d_K}(w_i+1)+S_{m,d_K}(2w_i)+w_i\leq 3S_{m,d_K}(2w_i)$.
It follows that $u_j'$ is a sum of $j$ terms which each have well-balanced constant $1,3S_{m,d_K}(2w_1),\dots,3S_{m,d_K}(2w_{j-1})$, respectively.
As the maximum of these terms is the last one, we can crudely estimate that 
\begin{align*}
w_j\leq S_{j,d_K}(3S_{m,d_K}(2w_{j-1}))\leq S_{r,d_K}(3S_{m,d_K}(2w_{j-1})).
\end{align*}
This is a first-order linear recurrence $w_j\leq \rho_1w_{j-1}+\rho_0$, with coefficients given by
\begin{align*}
\rho_1&:=\rho_1(m,r,d_K):=\begin{cases}
6mr&d_K=1\\
6
\left[d_K^{m-1}+d_K\frac{d_K^{m-1}-1}{d_K-1}\right]
\left[d_K^{r-1}+d_K\frac{d_K^{r-1}-1}{d_K-1}\right]&d_K>1,
\end{cases}\\
\rho_0&:=\rho_0(m,r,d_K):=\begin{cases}
(3m-2)r+1&d_K=1\\
\left[3d_K\frac{d_K^{m-1}-1}{d_K-1}\right]
\left[d_K^{r-1}+d_K\frac{d_K^{r-1}-1}{d_K-1}\right]
+d_K\frac{d_K^{r-1}-1}{d_K-1}&d_K>1.
\end{cases}
\end{align*}
The recurrence has the following bound,
\[w_j\leq \rho_1^{j-1}+\rho_0\frac{\rho_1^{j-1}-1}{\rho_1-1}.\]

From the Gram-Schmidt procedure, each $u_j'$ can be written as a linear combination of $u_1,\dots,u_{m'}$, and a suitable well-balanced constant for the scalars is $\B_j:=\B_j(m,r,d_K):=S_{j,d_K}(S_{m,d_K}(2w_{j-1}))$.
So we can write $u_j'=\sum_{i=1}^rc_{ij}b_i$, where each $c_{ij}$ is $\B_j$-well-balanced.
By (\ref{eq:distcomp}), we can estimate
\begin{align*}
\dist(v,H)
&\leq\sum_{i=1}^r\sum_{j=1}^{m'}\rabs{\frac{c_{ij}\lb v,b_i\rb}{\lb u_j',u_j'\rb}u_j'}
\leq\paren{\max_i\rabs{\lb v,b_i\rb}}
\sum_{i=1}^r\sum_{j=1}^{m'}\rabs{\frac{c_{ij}}{\lb u_j',u_j'\rb}u_j'}.
\end{align*}
Since $c_{ij},\lb u_j',u_j'\rb$ and $u_j'$ have well-balanced constants $\B_{j,d_K},S_{m,d_K}(2w_j)$ and $w_j$, respectively, the following constant,
\begin{align}\label{eq:theta1}
\T_1:=\T_1(m,r,d_K):=S_{rm,d_K}(\B_m+S_{m,d_K}(2w_m)+w_m),
\end{align}
is a well-balanced constant for the double sum.
\epf

\begin{remark*}\label{remark:technicalannoyance}
Fixing $m$ and $r$, it is clear that $\T_1(m,r,d_K)$ is an increasing function of $d_K$.
Accordingly, when we apply Lemma \ref{lem:peter5a}, it is sufficient (and convenient) to use an upper bound on $d_K$ as the third argument of $\T_1$.
\end{remark*}

\begin{lemma}\label{lem:peter5b}
Assume the hypotheses of Lemma \ref{lem:peter5a},  and additionally suppose $H\neq 0$ and $K=\Q$.
Then, there is an $\R$-basis $\set{v_i}$ of $H$, with entries in $\Z$, such that $\n{v_i}\leq A^{\T_2}$, where
\begin{equation}\label{eq:theta2}
\T_2:=\T_2(m):=m\cdot S_{m-1,1}\big(S_{(m-1)!,1}(m-1)+S_{(m-2)!,1}(m-2)+1\big).
\end{equation}
\end{lemma}
\bpf
Say $H=\mathrm{span}\set{b_1,\dots,b_{m'}}^\top$ with $b_1,\dots,b_{m'}$ linearly independent, so $\dim H=m-m'$.
Let $B\in\Mat_{m'\times m}(\cO_K)$ such that its $i$'th row is $b_i$.
Since this matrix has full rank, we can change the coordinates to some $C=(C_1|C_2)$ with $C_1\in\Mat_{m'\times m'}$ invertible.
Now, any $y\in H$ can be decomposed as $y=(y_1,y_2)\in\R^{m'}\times\R^{m-m'}$ with $y_1=-C_1\inv C_2y_2$.
It is straightforward to compute that $\det C_1$ has well-balanced constant $S_{m'!,1}(m')$, so $C_1\inv$ has well-balanced constant $S_{m'!,1}(m')+S_{(m'-1)!,1}(m'-1)$.
It follows that $-C_1\inv C_2$ has well-balanced constant $\A(m'):=S_{m',1}(S_{m'!,1}(m')+S_{(m'-1)!,1}(m'-1)+1)$.
Next, letting $y_2$ range through the standard basis vectors of $\R^{m-m'}$ yields a basis of elements $y\in K^m$ with well-balanced constant $\A(m')$.
Multiplying by the denominators of the first $m'$ entries yields a basis of integral vectors $y$ with well-balanced constant $(m'+1)\A(m')$.
The bound now follows from $m-m'>0$.
\epf

Denote by $\Sym_n$ the vector space of $n\times n$ symmetric matrices, and $\Pos_n$ the subspace of positive-definite matrices.
Fix non-empty open bounded sets $\cM,\cM^*$ such that $\cM^*\sub\cl{\cM^*}\sub\cM\sub\cl\cM\sub\Pos_n$, where the bar denotes topological closure.
For a matrix $Q\in\Mat_n(\Q)$ we denote by $\den(Q)$ the smallest positive integer $r$ such that $rQ\in\Mat_n(\Z)$, and we also define $\tilde Q:=\den(Q)\cdot Q=(\tilde Q_{ij})$.
If $Q$ is symmetric and positive-definite, 
then we let $\cQ:=\{\tilde Q_{jj}:1\leq j\leq n\}$ be the diagonal entries of $\tilde Q$, and 
$\cD:=\{\tilde Q_{ii}\tilde Q_{jj}-\tilde Q_{ij}^2:1\leq i<j\leq n\}$ be the $2\times 2$ diagonal determinants.
We say that a prime $p$ is $Q$-\ita{good} if $p$ is coprime to all elements in $\cQ$ and $-d$ is a quadratic non-residue modulo $p$ for each $d\in\cD$.

The next lemma will allow us to exchange the matrix $Q$ in (\ref{eq:blackbox}) with one that has better diophantine properties.

\begin{lemma}\label{lem:peter8}
Set $\T_3:=\T_3(n)$ and $\T_4:=\T_4(n)$ as in (\ref{eq:theta3}) and (\ref{eq:theta4}), respectively.
Let $L>c_0(\cM^*,\cM,n)$ (defined below), $D\geq 1$, and $M\geq \T_3D$.
Define $I:=[L,2L^D]$ and suppose $\cP\sub\set{(p^\nu,q^\nu):p,q\in I,1\leq\nu\leq n}$.
Let $Q\in\cM^*$.
Then there exists a nonzero subspace $H\sub\Sym_n$ (depending on $Q,\cP,$ and $M$) defined in (\ref{eq:Hdefn}) below, such that for every matrix $Q'\in H\cap\cM$, the following inclusion holds for all $(p^\nu,q^\nu)\in\cP$,
\begin{equation}\label{eq:theinclusion}
\cS(Q,q^\nu,p^\nu,M)\sub \cS(Q',q^\nu,p^\nu,\infty).
\end{equation}
Moreover, there exists a subset $\cP'\sub\cP$ with $\rabs{\cP'}\leq n(n+1)/2$ such that, setting
\begin{equation}\label{eq:Kdefn}
	K:=\Q((qp^{n-1})^{2\nu/n}:(p^\nu,q^\nu)\in\cP'),
\end{equation}
there exists $Q'\in H\cap\cM\cap\Mat_n(K)$; and if $K=\Q$, then we can find such a $Q'$ satisfying
\begin{equation}\label{eq:denQ'bound}
	\den(Q')\ll_{n,\cM,\cM^*} L^{\T_4D}.
\end{equation}
\end{lemma}
\bpf
For $\ga\in\Mat_n(\Z)$ and $m\in\N$, we define the following linear map,
\[B_{\ga,m}:\Sym_n\to\Sym_n:Q\mapsto\ga^\top Q\ga-m^{1/n}Q.\]
If we set
\begin{equation}\label{eq:Hdefn}
H:=\bigcap_{\substack{(p^\nu,q^\nu)\in\cP\\\ga\in \cS(Q,q^\nu,p^\nu,M)}}\ker B_{\ga,(qp^{n-1})^{2\nu}},
\end{equation}
then (\ref{eq:theinclusion}) is satisfied by construction.
Now, to each $B_{\ga,(qp^{n-1})^{2\nu}}$ we associate a matrix in $\Mat_{n(n+1)/2}$ which represents this map with respect to the coordinates of the standard basis of $\Sym_n$.
We write this basis as $\set{J_{ij}:1\leq i\leq j\leq n}$, where the $(i,j)$ and $(j,i)$ entries of $J_{ij}$ are $1$, and all other entries are zero.
Take a minimal set of rows $b_1^\top,\dots,b_r^\top\in\R^{n(n+1)/2}$, $r\leq n(n+1)/2$, of these matrices that generate $H^\perp$.
Let $\cP'$ be the set of corresponding pairs $(p^\nu,q^\nu)$ from (\ref{eq:Hdefn}), and define $K$ as in (\ref{eq:Kdefn}).

The $b_j$ have entries that are either in $\Z$, or else of the form $a-(qp^{n-1})^{2\nu/n}$, with $(p^\nu,q^\nu)\in\cP'$ and $a\in\Z$ satisfying $\rabs a\leq 2\max\ga_{ij}^2$; here, $\ga\in\cS(Q,q^\nu,p^\nu,M)$ is the matrix corresponding to the vector $b_j$ under consideration.
Since $Q\in\cl\cM$, a compact subset of $\Pos_n$, we have $Q=P^\top RP$, where $P\in O_n(\R)$ and $R$ is a diagonal matrix with eigenvalues $d_i$ satisfying $1\ll_\cM d_i\ll_\cM 1$.
By the bound in (\ref{eq:Sdefn}), this implies $(P\ga)^\top R(P\ga)\ll_\cM L^{2nD}$.
Defining $\tilde\ga:=P\ga$, this means $\tilde\ga_{ij}^2\ll_\cM L^{2nD}$.
Since $-1\leq P\leq 1$, we have $\ga_{ij}\ll_\cM L^{nD}$, so $a\ll_\cM L^{2nD}$.
We also clearly have $(qp^{n-1})^{2\nu/n}\ll_{n} L^{2nD}$.
We can estimate $\Gal(\cl K/\Q)\leq n^2(n^2-1)$, since $K$ is contained in the number field obtained by adjoining to $\Q$ at most $n(n+1)$ prime roots, and $\cl{\Q(p^{1/n})}=\Q(p^{1/n},\zeta_n)$ has degree at most $n(n-1)$. 
So by Lemma \ref{lem:Galois_bounding}, we have that $a-(qp^{n-1})^{2\nu/n}$ satisfies (\ref{eq:WBCdefn}) with $A\ll_{n,\cM} L^{6n^3(n^2-1)D}$.
Crucially, this holds because the Galois conjugates of $a$ are also $a$, as well as that the Galois conjugates of $p^{1/n}$ have absolute value $p^{1/n}$.

So by Lemma \ref{lem:peter5a} we have 
\[\dist(Q,H)\ll_{\cM} L^{6n^3(n^2-1)\T_1D}\max_j\rabs{\lb Q,b_j\rb},\]
where $\T_1=\T_1(n(n+1)/2,n(n+1)/2,n^2(n^2-1))$ is defined in (\ref{eq:theta1}); see the Remark after Lemma \ref{lem:peter5a} regarding the third argument.
We have $B_{\ga,(qp^{n-1})^{2\nu}}(Q)\ll_\cM L^{(2-M)D}$ by (\ref{eq:Sdefn}), so of course $\max_j\rabs{\lb Q,b_j\rb}$ satisfies the same bound, hence $\dist(Q,H)\ll_{\cM} L^{6n^3(n^2-1)\T_1D+(2-M)D}$.
This implies $\dist(Q,H)\ll_\cM L^{(6n^3(n^2-1)\T_1+2)D-M}$. 
Then, the following choice of $\T_3$,
\begin{equation}\label{eq:theta3}
\T_3:=\T_3(n):=6n^3(n^2-1)\T_1\big({n(n+1)}/2,{n(n+1)}/2,n^2(n^2-1)\big)+4,
\end{equation}
forces $(6n^3(n^2-1)\T_1+2)D<M-1$, which implies $\dist(Q,H)\ll_\cM L\inv$.
(A technical note: we could have simply required $6n^3(n^2-1)\T_1D+(2-M)D\leq -1$, for which it would suffice to impose $M\geq\T_3'(n)$ for some $\T_3'(n)$, say, rather than $M\geq\T_3D$ as is our current hypothesis. 
This alteration would slightly increase the numerical value of $\de_n$, but we opt to present the computation as above so that our presentation more closely mirrors \cite[Lemma 8]{BlomerMaga}).
Defining $d:=\dist(\cl{\cM^*},\cM^c)\gg_{\cM^*,\cM}1$, it is clear that $\dist(Q,H)\leq d/2$ implies $H$ intersects $\cM$ in some parallelepiped $E_\la$, where $\la\gg_{\cM^*,\cM} 1$ is the minimal length of any side of this polytope.
Hence, we will only consider sufficiently large $L>c_0(\cM^*,\cM,n)$.
By density of $\Q\sub K$, we have $H\cap\cM\cap\Mat_n(K)\neq\emptyset$.

Now, let us assume $K=\Q$.
The previous paragraph implies $H=\set0$ is impossible, so by Lemma \ref{lem:peter5b} there is an $\R$-basis $\set{v_i}$ of $H$, with entries in $\Z$, satisfying $\n{v_i}\ll_{\cM} L^{6n^3(n^2-1)\T_2D}$, where $\T_2=\T_2(n(n+1)/2)$ is defined in (\ref{eq:theta2}).
Consider the lattice $L:=\mathrm{span}_\Z\set{v_1,\dots,v_t}$, where $t=\dim H=n(n+1)/2-1$.
The projection of $E_\la$ onto each $v_i$ has some positive width $d_i\gg_{\cM,\cM^*}1$.
Let $b_i\in\Zpos$ satisfy $\n{v_i/b_i}<d_i$.
Such a $b_i$ which is minimally chosen satisfies $b_i\ll_{\cM,\cM^*} L^{6n^3(n^2-1)\T_2D}$.
Then there exists some $a_i\in\Z$ so that $\sum_{i=1}^ta_iv_i/b_i$ is in $E_\la$.
We've therefore constructed a $Q'\in H\cap\cM\cap\Mat_n(\Q)$ which satisfies (\ref{eq:denQ'bound}), with
\begin{equation}\label{eq:theta4}
\T_4:=\T_4(n):=\frac{n(n+1)}{2}6n^3(n^2-1)\T_2(n(n+1)/2).
\end{equation}
This completes the argument.
\epf

We will apply this next lemma to construct a suitable set of primes $\cP$ for use in (\ref{eq:blackbox}).


\begin{lemma}\label{lem:primesinintervals}
For every $\e\in(0,1/10)$, there exists $q_0(\e)$ so that for all $q>q_0(\e)$, if $(a,q)=1$, then there exists $t\in[q^{4.99},q^5]$ and some dyadic interval $[t/2^{i+1},t/2^i]\sub[t^{1-2\e},t]$ which contains $\gg_\e (t/2^{i+1})^{1/2}$ primes $p\equiv a\Mod q$.
\end{lemma}
\bpf
Define $\pi(x_0,x;q,a):=\#\set{p\in[x_0,x]:p\equiv a\Mod q}$.
By Xylouris \cite[Lemma 6.2 b)]{Xylouris}, for $q>q_0(\e)$ there exists $t\in[q^{4.99},q^5]$ such that $\pi(1,t;q,a)\geq t^{1-\e}/(\phi (q)\log t)$, where $\phi$ is the Euler totient function.
Then by Brun-Titchmarsh \cite{MontgomeryVaughan}, we have
\[\pi(1,t^{1-2\e};q,a)\leq \frac{2t^{1-2\e}}{\phi(q)\log(t^{1-2\e}/q)}.\]
It follows that
\[\pi(t^{1-2\e},t;q,a)
\geq\frac{t^{1-\e}}{\phi(q)}
\left[\frac{1}{\log t}-\frac{2t^{-\e}}{\log(t^{1-2\e}/q)}\right].\]
Since $t\geq q^{4.99}$ we have $-1/\log(t^{1-2\e}/q)\geq-1/\log(t^{1-(4.99)\inv-2\e})$, hence the bracketed quantity is at most $c/\log t$ provided $c<1-2t^{-\e}/(1-(4.99)\inv-2\e)$, which holds for sufficiently large $t$.
Thus,
\begin{equation*}
\pi(t^{1-2\e},t;q,a)\gg \frac{t^{1-\e}}{\phi(q)\log(t)}.
\end{equation*}

Next, we decompose $[t^{1-2\e},t]$ into dyadic intervals, $[t/2,t]\cup[t/4,t/2]\cup\cdots$ until the left endpoint arrives below $t^{1-2\e}$.
Because there are at most $\log t$ such intervals, it follows that for some $i$ we have 
\begin{align*}
\pi(t/2^{i+1},t/2^i;q,a)
&\gg \frac{t^{1-\e}}{\phi(q)\log^2(t)}
\gg t^{1/2},
\end{align*}
since $\phi(q)\log^2(t)\leq t^{1/4}\log^2(t)\ll t^{1/2-\e}$.
\epf

\section{A doubly recursive argument}\label{sec:recursion}

Here we utilize a doubly recursive argument to achieve a good bound on $\#\cS(Q,q^\nu,p^\nu,M)$ for suitable primes in suitable intervals.
This proposition concludes our diophantine investigations, and is based on \cite[Proposition 1]{BlomerMaga}.
Our addition is the incorporation of a Linnik-type theorem from Xylouris \cite{Xylouris} which improves the corresponding result in \cite{BlomerMaga} by providing explicitly computed constants.

\begin{proposition}\label{prop:recursion}
Let $L>c_1(n,\cM,\cM^*)$ (given below) and let $M,D_1\geq 1$ be fixed parameters satisfying (\ref{eq:Mbound}) and (\ref{eq:D1bound}).
Let $Q\in\cM^*$.
Then there exists $\cL$ satisfying
\begin{equation}\label{eq:ultimateLbound}
L^{D_1}
\leq\cL
\ll_{n,\cM,\cM^*} 
L^{\left[D_1{n\choose 2}10\T_4(n)\right]^{{n+1\choose 2}}},
\end{equation}
as well as a set of primes $\cP\sub[\cL,2\cL]$ satisfying $|\cP|\gg_{n,\cM,\cM^*}\cL^{1/2}$, such that
\begin{equation}\label{eq:ultimateSbound}
\#\cS(Q,q^\nu,p^\nu,M)\ll_{n,\cM,\cM^*} p^{\nu(n-2+\e)+\paren{1-\frac{1}{n}}\frac{1}{9.97}}
\end{equation}
for all $p,q\in\cP$ and $1\leq\nu\leq n$.
\end{proposition}
\bpf
For $0\leq j\leq n(n+1)/2$ we define 
\[I_j:=[L,2L^{D_1^{j+1}}],\qquad 
\cP_j:=\set{(p^\nu,q^\nu):p,q\in I_j,1\leq\nu\leq n},\]
and with this choice of $\cP_j$ let $H_j\sub\Sym_n$ be as in (\ref{eq:Hdefn}).
Attached to these data is a field $K_j$ and a matrix $Q_j\in\cM\cap\Mat_n(K_j)\cap H_j$ as in Lemma \ref{lem:peter8}.
We have $\Sym_n\supseteq H_0\supseteq H_1\supseteq\dots$.
Therefore we must have $H_i=H_{i+1}$ for some $i<n(n+1)/2$.
Since $Q_i\in H_i=H_{i+1}$, we can apply Lemma \ref{lem:peter8} with the parameters $\cP_{i+1},D_1^{n(n+1)/2+1}$, and
\begin{equation}\label{eq:Mbound}
M\geq \T_3D_1^{n(n+1)/2+1},
\end{equation}
where $\T_3:=\T_3(n)$ is defined in (\ref{eq:theta3}), and conclude by (\ref{eq:theinclusion}) that, for all $(p^\nu,q^\nu)\in\cP_{i+1}$, we have $\cS(Q,q^\nu,p^\nu,M)\sub\cS(Q_i,q^\nu,p^\nu,\infty)$.
By \cite[(6.2)]{BlomerMaga} this implies the following bound,
\begin{equation}\label{eq:recursionbound1}
\rabs{\cS(Q,q^\nu,p^\nu,M)}
\leq\rabs{\cS(Q_i,q^\nu,p^\nu,\infty)}
=0,\qquad p\neq q\in I_{i+1}\sm I_i,\frac{2\nu}n\notin\N.
\end{equation}

The remaining cases to consider are (i) $q=p$, and (ii) $q\neq p$, but $2\nu/n\in\N$.
The union of these cases is equivalent to $(qp^{n-1})^{2\nu/n}\in\N$.
Let $\cL_0:=L^{D_1^{i+1}}$, and define the interval $I_0^*:=[\cL_0,2\cL_0]$, as well as the following set of pairs of prime powers,
\[\cP_0^*:=\set{(p^\nu,q^\nu):p,q\in I_0^*,1\leq\nu\leq n,(qp^{n-1})^{2\nu/n}\in\N}.\]
We then apply Lemma \ref{lem:peter8} with the parameters $Q$, $D=1$, $M$ as in (\ref{eq:Mbound}), and $\cP_0^*$, which yields $H_0^*$ as in (\ref{eq:Hdefn}), and since $K=\Q$ in this case there exists some matrix $Q_0^*\in H\cap\cM\cap\Mat_n(\Q)$ which satisfies $\den(Q_0^*)\ll_{\cM,\cM^*}\cL_0^{\T_4}$, where $\T_4:=\T_4(n)$ is defined in (\ref{eq:theta4}).

Next, we suppose $0<j\leq n(n+1)/2$.
We will inductively construct $\cL_j$ and $Q_j^*\in\Mat_n(\Q)$ so that
\begin{equation}\label{eq:inductiveLjbound}
\cL_{j-1}^{{n\choose 2}9.98(1-2\e)\T_4}
\ll_{\cM,\cM^*}
\cL_{j}
\ll_{\cM,\cM^*}
\cL_{j-1}^{{n\choose 2}10\T_4}
\end{equation}
and
\begin{equation}\label{eq:inductiveLrequirement}
\#\set{Q_{j-1}^*\text{-good primes}\ p\in[\cL_{j},2\cL_{j}]}
\gg_{\cM,\cM^*}
\cL_{j}^{1/2}
\end{equation}
and
\begin{equation}\label{eq:inductiveQrequirement}
\den (Q_j^*)\ll_{\cM,\cM^*}\cL_j^{\T_4}
\end{equation}
hold; recall that here, we allow ourselves to take $\cL_j$ sufficiently large to guarantee the relative asymptotic growth.

We first construct $\cL_j$ which satisfies (\ref{eq:inductiveLjbound}) and (\ref{eq:inductiveLrequirement}).
Towards this, associated to the rational matrix $Q_{j-1}^*$ we have the sets $\cD_{j-1}^*$ and $\cQ_{j-1}^*$, which were defined immediately preceding Lemma \ref{lem:peter8}.
For a prime $p$ to be $Q_{j-1}^*$-good, we first require $\leg{-d}p=-1$ for all $d\in\cD_{j-1}^*$.
We construct a system of congruences which suffices to imply this.
We first impose $p\equiv -1\Mod 4$, which ensures $\leg{-1}p=-1$.
Now, list all the possible prime factors of any of the $d$'s.
Call such a prime $r$, and now we run through them, imposing a few conditions on $p$: 
if $r=2$, then impose $p\equiv-1\Mod 8$; 
if $r\equiv 1\Mod 4$, then impose $p\equiv 1\Mod r$; 
and if $r\equiv -1\Mod 4$, then impose $p\equiv -1\Mod r$.
By multiplicativity of the Legendre symbol and quadratic reciprocity, these constraints imply that each $\leg{r}p=1$, so $\leg{-d}p=-1$.
Note that by compactness, we have $\cD_{j-1}^*\sub[1,O(\cL_{j-1}^{2\T_4})]$ since $\den(Q_j^*)\ll_{\cM,\cM^*}\cL_{j-1}^{\T_4}$ by (\ref{eq:inductiveQrequirement}) for $j-1$; additionally, we have $\#\cD_{j-1}^*\leq{n\choose 2}$.
Hence, in order to satisfy this system of congruences, by the Chinese remainder theorem it suffices to satisfy a single congruence $p\equiv a_{j-1}\Mod{q_{j-1}}$ for some $q_{j-1}\asymp\cL_{j-1}^{{n\choose 2}2\T_4}$.
Now by Lemma \ref{lem:primesinintervals}, there exists $t_{j-1}\in[q_{j-1}^{4.99},q_{j-1}^5]$ and some dyadic interval $[t_{j-1}/2^{i_{j-1}+1},t_{j-1}/2^{i_{j-1}}]\sub[t_{j-1}^{1-2\e},t_{j-1}]$ which contains $\gg_\e (t_{j-1}/2^{i_{j-1}+1})^{1/2}$ primes $p\equiv a_{j-1}\Mod{q_{j-1}}$.
Finally, we choose $\cL_{j}:=t_{j-1}/2^{i_{j-1}+1}$, so 
(\ref{eq:inductiveLjbound}) holds.

For $Q_{j-1}^*$-goodness we also require $(\tilde Q_{\ell\ell},p)=1$ for each $\tilde Q_{\ell\ell}\in\cQ_{j-1}^*$.
As before, by compactness we have that $\cQ_{j-1}^*\sub[1,O(\cL_{j-1}^{\T_4})]$; additionally, we have $\#\cQ_{j-1}^*\leq n$.
Hence, the quantity of prime divisors of $\cQ_{j-1}^*$ is $\ll_{\cM,\cM^*}\T_4\log\cL_{j-1}$.
If we need to remove this many primes from $[\cL_{j},2\cL_{j}]$, then for $\cL_j$ sufficiently large, there would still remain $\gg_{\cM,\cM^*}\cL_{j}^{1/2}$ primes in $[\cL_{j},2\cL_{j}]$ which are $Q_{j-1}^*$-good; indeed, by (\ref{eq:inductiveLjbound}) we can estimate that 
\begin{align*}
\cL_{j}^{1/2}-\T_4\log\cL_{j-1}
\gg\cL_{j}^{1/2}-\left[{n\choose 2}9.98(1-2\e)\right]\inv\log\cL_{j}
\gg\cL_{j}^{1/2}.
\end{align*}
Thus $(\ref{eq:inductiveLrequirement})$ holds as well.

To finish the induction, we now construct $Q_j^*$ which satisfies (\ref{eq:inductiveQrequirement}).
Towards this we define the intervals $I_j^*:=[\cL_0,2\cL_j]$ and $\tilde I_j^*:=[\cL_j,2\cL_j]$, as well as the following sets of pairs of prime powers,
\begin{align*}
\tilde\cP_j^*&:=\set{(p^\nu,q^\nu):p,q\in \tilde I_j^*,1\leq\nu\leq n,(qp^{n-1})^{2\nu/n}\in\N,p,q\ \text{are}\ Q_{j-1}^*\text{-good}}\\
\cP_j^*&:=\cP_{j-1}^*\cup\tilde\cP_j^*.
\end{align*}
We take $H_j^*:=H_j(\cP_j^*,M)$ as in (\ref{eq:Hdefn}), where $M$ is as in (\ref{eq:Mbound}).
We then apply Lemma \ref{lem:peter8} with $Q$ and $P=\cP_j^*$, which yields a matrix $Q_j^*\in H_j^*\cap\cM\cap\Mat_n(\Q)$ which satisfies (\ref{eq:inductiveQrequirement}).
(Note that in the present case, the number field (\ref{eq:Kdefn}) is always $\Q$.)
This completes the induction.

We claim that $I_0^*\sub I_1^*\sub\cdots\sub I_{n(n+1)/2}^*\sub I_{i+1}\sm I_i$.
The inclusion $I_j^*\sub I_{j+1}^*$ is equivalent to $\cL_j\leq\cL_{j+1}$.
In order for $I_{n(n+1)/2}^*\sub I_{i+1}\sm I_i$, we must have $\cL_{n(n+1)/2}\leq L^{D_1^{i+2}}.$
We can estimate $\cL_{n(n+1)/2}
\ll_\W L^{
D_1^{i+1}
\left[{n\choose 2}10\T_4\right]^{n+1\choose 2}
},$
so if we choose
\begin{equation}\label{eq:D1bound}
D_1>\T_5:=\T_5(n):=\left[{n\choose 2}10\T_4\right]^{n+1\choose 2},
\end{equation}
then it is clear that
\begin{align*}
\cL_{n(n+1)/2}
\ll_{\cM,\cM^*} L^{D_1^{i+2}}L^{D_1^{i+1}(\T_5-D_1)}
\leq L^{D_1^{i+2}}.
\end{align*}
The factor $L^{D_1^{i+1}(\T_5-D_1)}$ kills the implied constant for sufficiently large $L$, so inequality holds for $L>c_1(n,\cM,\cM^*)$.

These interval inclusions imply $\Sym_n\supseteq H_0^*\supseteq H_1^*\supseteq\dots$, so we must have $H_k^*=H_{k+1}^*$ for some $0\leq k<n(n+1)/2$.
Since $Q_k^*\in H_k^*=H_{k+1}^*$, it follows from (\ref{eq:theinclusion}) that $\cS(Q,q^\nu,p^\nu,M)\sub\cS(Q_k^*,q^\nu,p^\nu,\infty)$ for all $(p^\nu,q^\nu)\in\tilde\cP_{k+1}^*$.
Since this set consists of powers of $Q_k^*$-good primes, we conclude from \cite[Lemma 7]{BlomerMaga} the following bound,
\begin{align}\label{eq:recursionbound2}
|\cS(Q,q^\nu,p^\nu,M)|
&\leq|\cS(Q_k^*,q^\nu,p^\nu,\infty)|
\ll_{\cM,\cM^*} p^{\nu(n-2+\e)},
\qquad (p^\nu,q^\nu)\in\tilde\cP_{k+1}^*,
p\neq q,
\end{align}
and by \cite[Lemma 6]{BlomerMaga},
\begin{align}\label{eq:recursionbound3}
|\cS(Q,p^\nu,p^\nu,M)|
&\leq|\cS(Q_k^*,p^\nu,p^\nu,\infty)|
\ll_{\cM,\cM^*} p^{\nu(n-2+\e)+\frac{1}{9.97}\paren{1-\frac{1}{n}}},
\qquad (p^\nu,p^\nu)\in\tilde\cP_{k+1}^*.
\end{align}
Finally, we choose $\cL:=\cL_{k+1}$ and $\cP:=\set{Q_k^*\text{-good primes}\ p\in[\cL_{k+1},2\cL_{k+1}]},$ so (\ref{eq:ultimateLbound}) holds.
And combining the estimates (\ref{eq:recursionbound1}), (\ref{eq:recursionbound2}) and (\ref{eq:recursionbound3}) implies (\ref{eq:ultimateSbound}).
\epf

\section{Proof of Theorem \ref{thm:1}}\label{sec:proofofthm1}

We apply Proposition \ref{prop:recursion} with the parameters
\[D_1=\T_5+1,\qquad M=\T_3(\T_5+1)^{{n+1\choose 2}+1},\qquad L=\la_F^{\frac{n(n-1)}4\eta},\]
where $\eta>0$ is some small constant to be specified in a moment.
This yields $\cL$ as in (\ref{eq:ultimateLbound}) and a corresponding prime set $\cP$ with $\rabs\cP\gg_\W\cL^{1/2}$.
Hence by (\ref{eq:blackbox}), 
\begin{equation}\label{eq:blackboxwithparticularparameters}
|F(g)|^2
\ll_\W
\la_F^{\frac{n(n-1)}4}
\paren{
    \cL^{-1/2}
    +\la_F^{-1/4}\cL^{n^3+M/2}
    +\cL^{-1+\paren{1-\frac{1}n}\frac{1}{9.97}+\e}
}.
\end{equation}
An quick computation reveals that the following constants,
\[a(n):=\frac{n(n-1)}4(\T_5+1),\qquad 
b(n):=\frac{n(n-1)}4\left[(\T_5+1){n\choose 2}10\T_4\right]^{n+1\choose 2},\]
satisfy
$\la_F^{\eta\cdot a(n)}\ll_\W\cL\ll_\W \la_F^{\eta\cdot b(n)},$
so (\ref{eq:blackboxwithparticularparameters}) becomes
\begin{equation*}
    |F(g)|^2
    \ll_\W
    \la_F^{\frac{n(n-1)}4}
    \paren{
         \la_F^{-\eta\xi_1(n)}
        +\la_F^{-\frac{1}4+\eta\xi_2(n)}
        +\la_F^{-\eta\xi_3(n,\e)}
    },
\end{equation*}
where $\xi_1(n):=a(n)/2$, and
\[\xi_2(n):=b(n)\paren{n^3+\frac{M}2},\qquad
\xi_3(n,\e):=a(n)\paren{1-\paren{1-\frac{1}n}\frac{1}{9.97}-\e}.
\]
If we choose $\eta=1/{4(\xi_2(n)+\xi_3(n,0))}$, then it follows that
\begin{equation}\label{eq:delta_ndefn}
\de_n:=\frac{\xi_3(n,0)}{8(\xi_2(n)+\xi_3(n,0))}
\end{equation}
is admissible in (\ref{eq:mainthmFbound}).

We now sketch the computation of the asymptotic lower bound $\de_n\gg n^{-cn^6}$.
Clearly we have $\de_n\gg \xi_3(n,0)/\xi_2(n)$.
Elementary calculations reveal the following estimates,
\[\xi_3(n,0)\gg \paren{\frac{n(n+1)}{2}!}^{c_1n^2},\qquad\xi_2(n)\ll n^{c_2n^4}\paren{\frac{n(n+1)}{2}!}^{c_3n^4},\]
for some positive absolute constants $c_i$.
The desired bound now follows from Stirling's approximation for $(n(n+1)/2)!$.

\section{Proof of Theorem \ref{thm:2}}\label{sec:proofofthm2}

We first provide two results which bound the solution sets of relevant quadratic forms.
These are Lemmas \ref{lem:3bBlomerPohl} and \ref{lem:53BM2}, which are explicit versions of \cite[Lemma 3(b)]{BlomerPohl} and \cite[Corollary 5.3]{BlomerMaga2}, respectively.
We then apply Lemma \ref{lem:53BM2} to bound $\#\cS(Q,q^\nu,p^\nu,M)$ in the case $n=3$, yielding $\de_3$.

We denote by $H(P)$ the height of a quadratic polynomial $P$, which is the maximum of the absolute values of the coefficients of $P$.

\begin{lemma}\label{lem:3bBlomerPohl}
For each $\de,D>0$ and each quadratic polynomial $P(x,y)\in\R[x,y]$ whose quadratic homogeneous part is positive definite with discriminant $\rabs\D\geq D$, the bound $\rabs{P(x,y)}\leq\de$ implies $\max(\rabs x,\rabs y)\ll_D(\de+1+H(P))^3$. 
\end{lemma}
\bpf
We write $P(x,y)=ax^2+bxy+cy^2+dx+ey+f$ and $\D=b^2-4ac<0$.
Without loss of generality, we assume $c\leq a$.
Arguing as in \cite[Lemma 3]{BlomerPohl}, we have
\[y^2\ll_D\frac{(\rabs{P(-\xi,-\eta)}+\de)4H(P)}{\D}+\eta^2,\]
where $\xi:=(be-2cd)/\D$ and $\eta:=(bd-2ae)/\D$.
Since $\eta,\xi\ll H(P)^2/\D$, the above estimate implies
\[y^2\ll_D H(P)^6+H(P)^4+H(P)^2+\de H(P).\]
This is at most $(\de+1+H(P))^6$, which implies the desired bound for $|y|$.

By \cite[(7.6)]{BlomerPohl}, we have
\[\rabs x
\leq |\xi|
+\frac{\rabs{b(y+\eta)}}{2\rabs a}
+\frac{1}{2\rabs a}\left[
4\rabs a
\paren{\de+\frac{\rabs\D(y+\eta)^2}{4\rabs a}+\rabs{P(-\xi,-\eta)}}
\right]^{1/2}.\]
Our assumption $c\leq a$ and $b^2-4ac=\D$ implies $a\geq\sqrt{-\D}/2$.
Using our bound on $\rabs y$, this implies
\begin{align*}
\rabs x
\ll_D
 H(P)^2
+H(P)^3+\Big[\de&+(1+\de+H(P))^6+H(P)^2(1+\de+H(P))^3\\
&+2H(P)^2(1+\de+H(P))^3+3H(P)^5+2H(P)^3+H(P)\Big]^{1/2},
\end{align*}
which is again $\ll_D(\de+1+H(P))^3$, as claimed.
\epf

\begin{lemma}\label{lem:53BM2}
Let $n\geq 2$.
Let $Q\in\Mat_n(\R)$ be a fixed symmetric positive definite matrix and let $X\geq 1$.
Let $0\leq k\leq n-2$, and let $x_1,\dots,x_k\in\Z^n$ be linearly independent of norm $\ll X$.
Let $q_0,\dots,q_k\in\R$ be bounded by $X^2$ and let $0<\de<X^{-N}$,
where $N:=N(k)> 73k+74$. 
Then,
\[\#\set{y\in\Z^n:y^\top Qy=q_0+O(X^2\de),x_j^\top Qy=q_j+O(X^2\de)\ \text{for}\ 1\leq j\leq k}
\ll X^{n-k-2+\e}.\]
\end{lemma}
\bpf
By \cite[Corollary 5.3]{BlomerMaga2}, the result holds for $N:=N(k)\geq k+2+14A(k+1)$, where the constant $A$ is inexplicitly provided by \cite[Corollary 4]{BlomerPohl}.
A straightforward computation reveals that $A>12C/7$ suffices, with the constant $C$ inexplicitly provided by \cite[Lemma 3(b)]{BlomerPohl}.
We computed in Lemma \ref{lem:3bBlomerPohl} that $C=3$ suffices.
\epf

Now, we will directly estimate $\#\cS(Q,q^\nu,p^\nu,M)$ using three applications of Lemma \ref{lem:53BM2}.
In Proposition \ref{prop:blackbox}, we choose $\cL=\la_F^\eta$, where $\eta>0$ is some constant which we will specify later, and let $\cP$ be the set of primes in $[\cL,2\cL]$.

For any $\ga\in\cS(Q,q^\nu,p^\nu,M)$, its first column $y_1\in\Z^n$ satisfies 
\[y_1^\top\paren{\frac{Q}{Q_{11}}}y_1=(qp^{n-1})^{2\nu/n}+O(\cL^{(2-M)\nu}).\]
Hence, we apply Lemma \ref{lem:53BM2} with the matrix $Q/Q_{11}$, as well as $X=(2\cL)^\nu$, $k=0$, $q_0=(qp^{n-1})^{2\nu/n}$, and $\de=2^{-2\nu}\cL^{-M\nu}$, where $M>N(0)$.
It follows that there are $\ll_\W\cL^{\nu+\e}$ possible choices for $y_1$.
Also, the second column $y_2\in\Z^n$ satisfies 
\[y_2^\top Qy_2
=(qp^{n-1})^{2\nu/n}Q_{22}+O(L^{(2-M)\nu}),
\quad y_1^\top Qy_2
=(qp^{n-1})^{2\nu/n}{Q_{12}}+O(L^{(2-M)\nu}).\]
So if we define $\kappa_2:=\max\set{Q_{11},\sgn(Q_{12})\cdot Q_{12}}$, then we can apply the Lemma with the matrix $Q/\kappa_2$, as well as $X=(2\cL)^\nu$, $k=1$, $q_0=(qp^{n-1})^{2\nu/n}Q_{22}/\kappa_2$, $q_1=(qp^{n-1})^{2\nu/n}Q_{12}/\kappa_2$, and again $\de=2^{-2\nu}\cL^{-M\nu}$, where this time we require $M>N(1)>147$.
Thus, there are $\ll_\W\cL^{\e}$ possible choices for $y_2$.
Similarly, we get that there are $\ll_\W\cL^\e$ possible choices for the third column of $\ga$.

We are now in a position to apply Proposition \ref{prop:blackbox}.
We argued that there are $\ll_\W \cL^{\nu+\e}$ different choices for $\ga$, provided $M>147$ in (\ref{eq:blackbox}).
By the prime number theorem we have $\rabs\cP\gg\cL^{1-\e}$, so by (\ref{eq:blackbox}) we get
\[\rabs{F(g)}^2\ll_\W \la_F^{\frac{n(n-1)}4}
\paren{
\la_F^{\eta(-1+\e)}
+\la_F^{-1/4}\la_F^{\eta(n^3+M/2)}
+\la_F^{\eta(-1+\e)}
}.\]
If we choose $\eta=1/(4+4n^3+2M),$
then it follows that
\[\de_3=\frac{1}{8+8n^3+4M}\]
is admissible in (\ref{eq:mainthmFbound}).

\section*{Acknowledgements}
The author would like to thank P\'eter Maga for his stellar mentorship and guidance, as well as useful discussions and feedback on previous versions of this paper.
In addition, he would like to thank the anonymous referee for the helpful suggestions regarding the exposition.
He would also like to express his gratitude towards the Budapest Semesters in Mathematics program for providing the framework under which this research was conducted.

\end{document}